%This is an amstex2.1-file.

%
% This file was prepared using
% AMS-TeX 2.1 with AMS-Fonts 2.1
% and the style file amsppt.sty version 2.1
%
\documentstyle{amsppt}
%\magnification=1100
%%\pagewidth{6.4in}\vsize8.5in\parindent=6mm
%%\parskip=3pt\baselineskip=14pt\tolerance=10000\hbadness=500
%%\NoRunningHeads
\loadbold
\topmatter
\title
Pointwise convergence of lacunary spherical means
\endtitle
\author
Andreas Seeger \ \ \  Terence Tao \ \ \ James Wright
\endauthor
%%\date  April 2, 2002\enddate
\leftheadtext{A. Seeger, T. Tao, J. Wright}
\abstract
We show 
that if $f$ is locally in $L\log\log L$ then the  lacunary spherical means
converge almost everywhere.
The argument given here is a model case for more general results 
 on singular maximal functions and Hilbert transforms along plane
 curves \cite{6}.
  \endabstract
\thanks
The first author is supported in part by a grant from the National Science
Foundation.  The second author is a Clay Prize fellow
and is supported by the Sloan and Packard foundations.
  \endthanks
\address
Department of Mathematics,
University of Wisconsin, Madison, WI 53706-1388, USA
\endaddress
\email seeger\@math.wisc.edu\endemail
\address
Department of Mathematics,
University of California, Los Angeles, CA  90095-1555, USA
\endaddress
\email tao\@math.ucla.edu\endemail
\address
Department of Mathematics and Statistics,
University of Edinburgh,
King's Buil\-dings,
Mayfield Rd.,
Edinburgh EH3 9JZ, U.K.
\endaddress
\email wright\@maths.ed.ac.uk\endemail

\endtopmatter
\document

\def\vth{\vartheta}

\def\Q{\Cal Q}

\def\emph#1{{\it #1 }}

\define\th{\Theta}

%\define\prd{{\text{\rm prod}}}

\define\dist{{\text{\rm dist}}}

\define\inn#1#2{\langle#1,#2\rangle}
\define\biginn#1#2{\big\langle#1,#2\big\rangle}

\define\lcontr{\rfloor}
\define\lco#1#2{{#1}\lcontr{#2}}
\define\lcoi#1#2{\imath({#1}){#2}}
\define\rco#1#2{{#1}\rcontr{#2}}

\define\bin#1#2{{\pmatrix {#1}\\{#2}\endpmatrix}}
\define\meas{{\text{\rm meas}}}

\define\lc{\lesssim}
\define\gc{\gtrsim}

%Greek letters

\define\ka{\kappa}
            
\define\la{\lambda}

\define\fQ{{\frak Q}}

%Attn:\fi can't be defined.

\define\bbR{{\Bbb R}}

\define\cA{{\Cal A}}

\define\cM{{\Cal M}}

\define\cQ{{\Cal Q}}

\define\cV{{\Cal V}}

%roman letters with a tilde

%roman letters with a bar
%\define\abar{{\bar a}} \define\Abar{{\bar A}}
%\define\bbar{{\bar b}} \define\Bbar{{\bar B}}
%\define\cbar{{\bar c}} \define\Cbar{{\bar C}}
%\define\dbar{{\bar d}} \define\Dbar{{\bar D}}
%\define\ebar{{\bar e}} \define\Ebar{{\bar E}}
%\define\fbar{{\bar f}} \define\Fbar{{\bar F}}
%\define\gbar{{\bar g}} \define\Gbar{{\bar G}}
%\define\hBar{{\bar h}} \define\Hbar{{\bar H}}
%\define\ibar{{\bar i}} \define\Ibar{{\bar I}}
%\define\jbar{{\bar j}} \define\Jbar{{\bar J}}
%\define\kbar{{\bar k}} \define\Kbar{{\bar K}}
%\define\lbar{{\bar l}} \define\Lbar{{\bar L}}
%\define\mbar{{\bar m}} \define\Mbar{{\bar M}}
%\define\nbar{{\bar n}} \define\Nbar{{\bar N}}
%\define\obar{{\bar o}} \define\Obar{{\bar O}}
%\define\pbar{{\bar p}} \define\Pbar{{\bar P}}
%\define\qbar{{\bar q}} \define\Qbar{{\bar Q}}
%\define\rbar{{\bar r}} \define\Rbar{{\bar R}}
%\define\sbar{{\bar s}} \define\Sbar{{\bar S}}
%\define\tbar{{\bar t}} \define\Tbar{{\bar T}}
%\define\ubar{{\bar u}} \define\Ubar{{\bar U}}
%\define\vbar{{\bar v}} \define\Vbar{{\bar V}}
%\define\wbar{{\bar w}} \define\Wbar{{\bar W}}
%\define\xbar{{\bar x}} \define\Xbar{{\bar X}}
%\define\ybar{{\bar y}} \define\Ybar{{\bar Y}}
%\define\zbar{{\bar z}}  \define\Zbar{{\bar Z}}

%roman letters with a hat

%%%% \input def

%%\head{\bf 1.Introduction}
%%\endhead
\subheading{1. Introduction}
Let $f$ be a locally integrable function on $\Bbb R^d$ where $d\ge 2$.
For any  integer $k$ let $\cA_k f(x)$ be
the spherical average of $f$ over the  sphere
of radius $2^k$ in $\Bbb R^d$   centered at $x$; i.e.
$$
\cA_k f(x)=\int f(x+2^{k}y) d\theta(y);
$$
here $d\theta$ denotes the normalized Lebesgue measure  on the unit sphere.
Clearly we have $\lim_{k\to -\infty} \cA_k f(x)=f(x)$ for all $x$ if $f$ is  continuous
everywhere.
Moreover by results  of C. Calder\'on \cite{1} and Coifman and Weiss \cite{4}
we have
$$\lim_{k\to-\infty} \cA_k f(x)=f(x) \quad\text{   almost everywhere }
\tag 1.1 $$ if $f$ is locally in $L^p$ for $p>1$.
It is well known \cite{7} that such results are equivalent with a weak type
$(p,p) $ bound for the local maximal function
$\cM f(x)=\sup_{k\le 0}|\cA_k f(x)|$
and the above mentioned authors showed that the maximal operator is bounded on $L^p$.
It is still unknown whether  (1.1) holds for $f\in L^1$ and, equivalently,
whether the maximal function is of weak type (1,1).
However we have 
\proclaim{Theorem}
%%For $t\in (0,\infty) $ let
%%$$\Phi(t)= t\log\log(10+t).
%%
%%$$
Let $\cM$ be the global lacunary spherical maximal operator defined by
 $$\cM f(x)=\sup_{k\in \Bbb Z}|\cA_k f(x)|.
\tag 1.2
$$
There is a constant $C$ so that for all measurable  functions $f$ 
and all $\alpha>0$ the inequality
$$
\big|\{x\in \Bbb R^n: \cM f(x)>\alpha\}\big|
\le
%%\int \Phi\big(\frac{C|f(x)|}{\alpha}\big) dx
\int \frac{C|f(x)|}{\alpha}\log\log\big(e^2+\frac{C|f(x)|}{\alpha}\big) dx
\tag 1.3
$$
holds.
\endproclaim

As a corollary one obtains the pointwise convergence result (1.1) if
%%$\lim_{k\to-\infty}\cA_kf(x)=f(x)$ holds almost everywhere
$f$  belongs locally to
$L\log \log L$.

The inequality  for the lacunary spherical maximal function is a special case of more general
 results in \cite{6} which apply  to 
operators such as maximal averages and Hilbert transforms 
along plane curves which are homogeneous with respect to some
family of a  nonisotropic dilation. The presence of the nonisotropic dilation structure causes additional difficulties and 
therefore it seems adequate
to present the less technical proof of the Theorem above  separately.
 The main idea from \cite{6}  is still present.

Concerning weak type inequalities for classes near $L^1$ we mention two previous results 
for the operator $\cM$. First
 Christ and Stein \cite{3}  showed (combining Calder\'on-Zygmund arguments with Yano's extrapolation theorem)
that  $\cM$ maps $L\log L(Q_0)$ (for a unit cube $Q_0$) to weak $L^1$. This result
applies to more general maximal operators associated to Borel measures, whose Fourier transform decays at infinity. 
Moreover Christ \cite{2}    showed the harder result
that $\cM$  maps the Hardy space $H^1(\bbR^d)$ to weak $L^1$
(for  weak type $(p,p)$ endpoint bounds for related maximal operators 
see also the recent paper \cite{5}).
The condition
$f\in H^1$ means that $f$
 has  some  rather substantial cancellation. Concerning size estimates 
note there is 
 a  restriction if  $f\in H^1$ is single signed in an open  ball; namely 
then $f$ belongs to $L\log L(K) $ for all  compact subsets 
$K$ of this ball
({\it cf.}\cite{8, \S I.5.2 (c)}).

{\it Notation.} 
For two quantities $A$ and $B$ we write 
 $A\lc B$ or $B\gc A$ if there exists an absolute 
 positive constant $C$ so that $A\le C B$. 
The Lebesgue measure of a set $E$ will  be denoted by  $|E|$.

\subheading{\bf 2. Length and thickness}
We say that a set $E$ is {\it granular} 
if $E$ is a finite  union of dyadic cubes.
%Let $q$ be a dyadic cube and let $E$ be a granular set. 
%The assumption of being granular simplifies various 
%stopping time arguments and is not a severe
%restriction since we observe
%
For a granular   set $E$   we define a dyadic version of a {\it one-dimensional Hausdorff content} 
or simply  {\it 'length' } $\la(E)$
  to be
$$ \la(E) := \inf_{\cQ} \sum_{Q \in \cQ} l(Q)\tag 2.1$$
where $\cQ$ ranges over all finite collections $\cQ$ of dyadic cubes  with
$E\subset\cup_{Q\in\cQ}Q$,
 and $l(Q)$ denotes the side-length of $Q$.
%%Note that this definition differs from the usual definition of a $1$ dimensional 
%% Hausdorff measure as $\la(E)\le l(Q)$ if $E$ is contained in the dyadic cube $Q$.

Next if $E$ is granular, we define the {\it thickness} $\th(E)$ to be
$$ \th(E) := \sup_{Q} \frac{|E \cap Q|}{l(Q)}\tag 2.2$$
where $Q$ ranges over all dyadic cubes.
Clearly, if $E$ is contained in a dyadic cube $q$  it is sufficient to let $Q$ in (2.2) range over all dyadic subcubes of $q$.

The quantities of length and thickness are complementary.  For instance, it is immediate from the definitions that one has
$$ |E| \leq \la(E) \th(E).
\tag 2.3$$

The bound (2.3) can be attained, for instance if $E$ is a dyadic box.
More generally if $C\ge 1$, we call $E$ a \emph{generalized box} with admissible deviation $C$
 if one has
$$\la(E) \th(E)\le C|E|.\tag 2.4$$

In the proof of the weak type $ L\log\log L$  inequality 
the quantity $\la(E)$ will control the size of the exceptional set, while $\th(E)$ 
will control the $L^2$ norm of the maximal function outside of the exceptional set.  
Inequalities such as (2.4) will be crucial  for balancing the two estimates. 
There is also an intermediate range of scales in which neither the of the quantities
 $\la(E)$ or $\th(E)$ quantities is favorable, and one will just use $L^1$ estimates 
for that portion.

Now unfortunately (2.4) is not always satisfied but the following 
proposition can be used to efficiently decompose a granular set  into generalized boxes of 
different lengths.

\proclaim{Proposition}  Let  $E$ be a granular set. Then
there exists a decomposition $E =F \cup G$ into disjoint granular subsets  $F$, $G$
such that
$$ \la(F) \leq \frac{1}{2} \la(E)\tag 2.5$$
and
$$ \th(G) \leq 8\frac{|G|}{\la(E)}.\tag 2.6$$
In particular, $G$ is a generalized box with admissible  deviation of at most $8$.
\endproclaim

\demo{\bf Proof}
Fix $E$
%; we shall write $\La$ for $\la(E)$.  
and define the {\it critical thickness} $\vth(E)$ 
to be the
largest non-negative number $r$ such that the inequality
$$  r\la(E) \leq 2r\sum_{Q \in \cQ} l(Q) +
\big|E \backslash \bigcup_{Q \in \cQ} Q\big| \tag 2.7$$
holds for all finite collections $\cQ$ of dyadic cubes (here the empty collection is admitted).
% contained in $q$.
Equivalently, one can define $\vth(E) $ by
$$ \vth(E)  := \inf_{\cQ} \frac{ |E \backslash \bigcup_{Q \in \cQ} Q| }{(\la(E) - 2\sum_{Q \in \cQ} l(Q))_+}.\tag 2.8$$
Clearly $\vth(E) \le |E|/\la(E)$. Observe also  that  $\vth(E)  > 0$.  This follows
 because $|E \backslash \bigcup_{Q \in \cQ} Q|$ is bounded away from
zero whenever $\sum_{Q \in \cQ} l(Q) \leq \la(E)/2$ (thanks to the
 hypothesis that $E$ is granular).
Moreover, again  since $E$ is granular, there exists a finite collection $\cQ_1$ of dyadic cubes such that
$$  \vth(E) \la(E) = 2\vth(E)  \sum_{Q \in \cQ_1} l(Q) + |E_*| \tag 2.9$$
where $E_*$ is the set
$$ E_* := E \backslash \bigcup_{Q \in \Q_1} Q.
\tag 2.10
$$
We claim that
$$ \th(E_*) \leq 2\vth(E) .\tag 2.11$$
Indeed, suppose that there existed a cube $Q'$ such that
$$ |E_* \cap Q'| >2 \vth(E)  l(Q').\tag 2.12$$
Then  $Q' \not \in \Q_1$.  If we apply (2.7) 
to the collection $\Q_1
\cup \{Q'\}$ we obtain
$$  \vth(E) \la(E) \leq 2\vth(E) (l(Q') + \sum_{Q \in \Q_1} l(Q)) + |E_*| - |E_* \cap Q'|,$$
but this contradicts (2.9) and (2.12).  This proves (2.11).

We now use a recursive construction to obtain large subsets of $E$ of bounded thickness.

\proclaim{Lemma} Let $r>0$. For any dyadic cube $I$, there exists a (possibly empty) collection
$\Q[I]$ of disjoint dyadic cubes in $I$ and a granular set $E[I] \subset E \cap I$ such that
$$\th(E[I]) \leq 2r \tag 2.13$$
and
$$ 2|E[I]| \geq 2r \sum_{Q \in \Q[I]} l(Q) + |(E \cap I) \backslash
\bigcup_{Q \in \Q[I]} Q|. \tag 2.14$$
\endproclaim

The Lemma is proved  by induction on the sidelength of $I$.  If $l(I) 
\leq (2r)^{\frac{1}{d-1}}$, the 
lemma follows simply by taking $E[I] = E \cap I$ and $\Q[I]$ to be empty.
Now 
%suppose $l(I) > (2r)^{d-1}$
fix a dyadic cube $I$ and suppose  that the lemma has been proven for all 
proper dyadic subcubes 
$I'$.  Partition $I$ into $2^d$ sub-cubes $I_1, \ldots, I_{2^d}$ of
side-length $\frac{1}{2} l(I)$.  By the inductive hypothesis,
we may construct collections $\Q[I_j]$ and
sets $E[I_j]$ for $j=1, \ldots, 2^d$
satisfying the properties of the lemma.
We distinguish two cases. 
First if $|\bigcup_{j=1}^{2^d} E[I_j]| \leq 2r l(I)$ then we 
define $E[I] := \bigcup_{j=1}^{2^d} E[I_j]$ and
$\Q[I] := \bigcup_{j=1}^{2^d} \Q[I_j]$.
Next if $|\bigcup_{j=1}^{2^d} E[I_j]| > 2r l(I)$ then we
take $E[I]$ to be a granular subset of
$\bigcup_{j=1}^{2^d} E[I_j]$ of measure at least $r l(I)$, and at most $2rl(I)$
 and take  $\Q[I]$ to be the singleton set $\{I\}$.  The properties (2.13/14) are not hard to check in both cases (see also \cite{6} for a more detailed description of a variant).
\enddemo

\demo{Proof of the Proposition, cont}
Since $E$ is granular there is a dyadic cube $q$ so that $E$ is contained in it.
We apply the Lemma with $I = q$ and $r=\vth(E)$. We thus  find a set $E[q]$ and a collection
$\Q[q]$ obeying the properties in the lemma.  We now set
$ G := E_* \cup E[q]$ 
and
$ F := E \backslash G.$

To show (2.5) we observe
$ F \subset E \backslash E_* \subset \bigcup_{Q \in \Q_1} Q,$
so
$ \la(F) \leq \sum_{Q \in \Q_1} l(Q).$
But by (2.9) this sum must be less than or equal to $\lambda(E)/2$,
which gives (2.5).

To show (2.6) we first observe that
$$ |G| \geq |E[q]|
\geq \frac{1}{2} (2\vth(E)) \sum_{Q \in \Q[q]} l(Q) + \Big|E \backslash \bigcup_{Q \in
\Q[q]} Q\Big|$$
by (2.14), since now $r=\vth(E)$.  By (2.7) we thus see that
$ |G| \geq \lambda(E) \vth(E)/2$.
Since
$ \th(G) \leq \th(E_*) + \th(E[q]) \leq 2\vth(E)  + 2\vth(E)  = 4\vth(E) ,$
we see that $\th(G)\le 8|G|/\la(E)$ which is (2.6).\qed

\enddemo

\subheading{\bf 3. Basic reductions}
We say that a function $f$ is granular if $f=\sum_\nu c_\nu \chi_{E_\nu}$
where $c_\nu\in \Bbb C$,  the sum is finite, the sets $E_\nu$ are granular and mutually disjoint.
In order to prove the Theorem  we may restrict ourselves  to  granular functions
since every measurable function which is finite almost everywhere is the almost
everywhere limit of a monotone sequence of granular functions.

Let now $f$ be a granular function and we shall estimate the size of
$\{x: \cM
 f(x)>\alpha\}$. We perform a
standard Calder\'on-Zygmund decomposition
at height $1$ for the function $\Phi(|f|/\alpha)$, where $\Phi(t)=t\log\log(e^2+t).$
This can be done via the Whitney decomposition theorem applied to the open set
$$\Omega= \{ x: M_{HL}(\Phi(|f|/\alpha))(x)>1\}$$
where $M_{HL}$ is the Hardy-Littlewood maximal function. We denote by
$\fQ$  the set of  Whitney
cubes arising in this fashion. By possibly subdividing each cube $Q$ into cubes of length $2^{-10}l(Q)$  we may assume
that 
$ 2^{8}l(Q)\le d^{-1/2}\dist(Q, \bbR^d\setminus \Omega_\alpha)\le 2^{10} l(Q)$.
From this it is easy to see that we may subdivide the family $\fQ$ 
into subcollections
$\fQ_1,\dots, \fQ_{N(d)}$ with the property that in each $\fQ_i$ the double cubes $Q^*$ are pairwise disjoint.
%\footnote{Alternatively use Terry's argument here.}

We shall slightly modify the definition of our maximal operator.
Let $d\sigma_k$ be the measure given by
$$\inn{f}{d\sigma_k}=\int f(2^{k }y)\chi(y) d\theta(y)$$
where $\chi$ is supported on a ball of radius $\le 1/2$ (so that the support of $d\sigma_k$
is contained in the sphere of radius $2^k$ centered at the origin, but  does not contain
antipodal  points on this sphere). We only need to consider
the
 maximal operator $M$ given by
$
M f(x)= \sup_{k\in \Bbb Z}
 |d\sigma_k* f(x)|.$

Now let
$g(x)=  f(x)$ if $|f(x)|\le  2^{10}\alpha$ and $g(x)=0$  otherwise.
Then $f-g$ is supported in $\Omega$.
Since $|g(x)|\lc\alpha$ the $L^2$ boundedness of the maximal operator and Chebyshev's inequality can be used to show that
$$
\Big|\{x: |M g(x) |> \alpha/2\}\Big|\lc
\int\frac{|f(x)|}{\alpha} dx \lc
\int\Phi(\frac{|f(x)|}{\alpha}) dx.
\tag 3.1
$$

Also if $\widetilde \Omega$ denotes the union of the tenfold expanded cubes then
$$|\widetilde \Omega|\lc|\Omega|\lc
\int\Phi(\frac{|f(x)|}{\alpha}) dx
\tag 3.2$$
and thus it suffices to show that
% measure of the set
%$\{x\notin \widetilde \Omega: |M[f-g](x)|>\alpha/2\}$ is controlled by 
%$\int\Phi(|f(x)|/\alpha) dx$.
%
$$
\Big|\{x\notin \widetilde \Omega: |M[f-g](x)|> \alpha/2\}\Big|\lc
\int\Phi(\frac{|f(x)|}{\alpha}) dx;
\tag 3.3
$$
(3.1), (3.2), (3.3) imply the assertion of the Theorem.
In order to prove (3.3) we split the function $f-g$ further.
For $n=10,11,\dots$ let
$$
E^n=\{x\in \Omega: 2^{n}\alpha<  |f|\le 2^{n+1}\alpha \}
$$
and let
$f^n(x)= f(x)\chi_{E^n}(x)$. Now $E^n=\cup_{q\in \fQ} E^n_q$ where $E^n_q=E^n\cap q$. Note that the sets $E^n_q$ are granular since
 $f$ was assumed to be a granular function.
We use Proposition 2.1 iteratively   to decompose $E^n_q$ further, namely 
$$
E^n_q=\bigcup_{\nu=0}^\infty E_q^{n,\nu}\cup F^n_q
$$where $F^n_q$ has measure zero and will be ignored in what follows,
 and where  each $E_q^{n,\nu}$ is a generalized box, the sets $E_q^{n,\nu}$ are mutually
disjoint, and $\lambda(E^{n}_q\setminus \cup_{\nu=1}^m E^{n,\nu}_q)\le 2^{-m}\lambda(E^{n}_q)$.
Thus also 
$$|E^{n,\nu}_q|\le
\la(E^{n,\nu}_q) l(q)^{d-1}\le
 2^{1-\nu}|q|.\tag 3.4$$
Finally define
$$f^{n,\nu}_q(x)= f(x) \chi_{ E^{n,\nu}_q}(x).\tag 3.5
$$

We shall first handle the terms with $\nu\ge n^2$ and show that
$$
\Big|\{x: \sup_k |\sum\Sb q\endSb\sum_{n}\sum_{\nu\ge n^2} 
f^{n,\nu}_q *d\sigma_k(x)|> \alpha/4\}\Big|
\lc
\int\Phi\big(\frac{|f(x)|}{\alpha}\big)dx. \tag 3.6
$$

 In view of (3.4) these  terms are again 
easy to handle  by an $L^2$ estimate.
By Chebyshev's inequality and  the $L^2$ boundedness of the  maximal operator we get
$$
\align
&\Big|\{x: \sup_k |\sum\Sb q\endSb\sum_{n}\sum_{\nu\ge n^2} 
f^{n,\nu}_q *d\sigma_k(x)|> \alpha/4\}\Big|
\lc \alpha^{-2}
\Big\|
\sum\Sb q\endSb\sum_{n}\sum_{\nu\ge n^2} 
f^{n,\nu}_q \Big\|_2^2
\\
&\lc \alpha^{-2}\int\Big |
\sum\Sb q\endSb\sum_{n}\sum_{\nu\ge n^2} 
2^n\alpha\chi_{E^{n,\nu}_q}(x) \Big|^2 dx
\lc 
\sum\Sb q\endSb\sum_{n}\sum_{\nu\ge n^2} 
2^{2n}|E^{n,\nu}_q|
\endalign$$
by the disjointness of the sets $E^{n,\nu}_q$. By (3.4) the last expression is bounded by
a constant times
$$
%%\align&
\sum\Sb q\endSb\sum_{n}\sum_{\nu\ge n^2} 
2^{2n} 2^{-\nu}|q|
\lc
\sum\Sb q\endSb|q|
%%\qquad\lc \big|\{x: M_{HL} [\Phi(|f|/\alpha)](x) >1\}\big|
 \lc \int \Phi\big(\frac{|f(x)|}{\alpha}\big) dx
%%\endalign
$$
which yields  (3.6).

We are left with the consideration of terms with $\nu<n^2$ in the complement of the 
set $\widetilde \Omega$.
Since for $2^k\le l(q)$ the convolution $d\sigma_k* f_q^{n,\nu}$ is supported in
$\widetilde \Omega$ we are reduced to verify that
$$
\Big|\{x: \sup_k |\sum\Sb q:\\2^k>l(q)\endSb\sum_{n}\sum_{\nu=1}^{n^2} f^{n,\nu}_q *d\sigma_k(x)|> \alpha/4\}\Big|\lc
\int\Phi(\frac{|f(x)|}{\alpha}) dx.
\tag 3.7
$$
This will be carried out in the next section. 
Clearly (3.6) and (3.7) imply the desired estimate (3.3).

\subheading{\bf 4. Proof of (3.7)}
For each $n,\nu, q$ let  $k^{n,\nu}_q$ be unique  integer for which
$$2^{k^{n,\nu}_q-1}<\max\big\{l(q), \big(2^n \log(10+n)\th(E^{n,\nu}_q)\big)^{1/(d-1)}\big\}
<2^{k^{n,\nu}_q}.
\tag 4.1
$$
We consider the contribution to the case $k\le k^{n,\nu}_q$. This contribution is supported
inside the set $\cV=\widetilde \Omega\cup\cV_1$ where
$$\cV_1= \bigcup_q \bigcup_n\bigcup_\nu \bigcup_{k: l(q) \leq 2^k \leq 2^{k^{n,\nu}_{q}}}
(E^{n,\nu}_q + S_k) \tag 4.2
$$
and $E^{n,\nu}_q + S_k$ is the Minkowski sum of the set $E^{n,\nu}_q$ and the sphere
$S_k$ of radius $2^k$.

By covering $E^{n,\nu}_q$ efficiently by cubes in $q$ ({\it cf.} the definition of $\la$)
 we see that the inner union has measure at most
$$ \lesssim \la(E^{n,\nu}_q) \th(E^{n,\nu}_q) 2^n\log(10+n)
\approx |E^{n,\nu}_q|2^n\log(1+n)
$$
since $E^{n,\nu}_q$ is a generalized box.
Thus
$$
|\cV| \lesssim|\widetilde \Omega|+ \sum_q \sum_n \sum_\nu|E^{n,\nu}_q|
2^n\log(10+n)
\lc  \int \Phi\big(\frac{|f(x)|}{\alpha}\big) dx
\tag 4.3
$$
by the disjointness and  definition of the sets $E^{n,\nu}_q$ .

Next let (for $n\ge 10$)  
$$\ka_n= 100 \log_2 n \tag 4.4 $$ and
 we consider the contribution of the scales
$$ k^{n,\nu}_q < k \leq k^{n,\nu}_q + \ka_n. \tag 4.5$$
For this case we replace the sup by the sum and
 use Chebyshev's inequality  in $L^1$ to estimate
$$
\align
&\meas\Big(\{x: \sup_k
\sum_{n,\nu}
\sum\Sb q:  k^{n,\nu}_q < k \\ \leq k^{n,\nu}_q + \ka_n \endSb
|f^{n,\nu}_q *d\sigma_k(x)|> \alpha/8\}\Big)
\\
&\le 8\alpha^{-1} \Big\|
 \sup_k
\sum_{n,\nu}
\sum\Sb q:  k^{n,\nu}_q < k \leq \\k^{n,\nu}_q+\ka_n \endSb
f^{n,\nu}_q *d\sigma_k\Big\|_1
%%\\&
\le 8\alpha^{-1}
\sum_q
\sum_{n,\nu}
\sum \Sb k:  k^{n,\nu}_q < k \leq \\k^{n,\nu}_q+\ka_n \endSb
\|f^{n,\nu}_q \|_1
\\
&\lc 
\sum_{n,\nu}  \log \log (10+2^n) \sum_q\int_q 2^n \chi_{E_q^{n,\nu}}(x) dx
%%\\&
\lc\int\Phi\big(\frac{|f(x)|}{\alpha}\big) dx.
\tag 4.6
\endalign
$$

For the remainder, we shall actually show that
$$ \meas\Big(\Big \{ x: \sup_k \big| \sum_q \sum\Sb n, 
\nu:\\ k > k_{q}^{n,\nu} + \ka_n
\\ \nu<n^2
\endSb
f^{n,\nu}_q*d\sigma_k (x) \big| >  \alpha/8 \Big\}\Big)
\lesssim \alpha^{-1}\int|f(x)| dx
\tag 4.7$$
and the right hand side is of course controlled by
$\int \Phi(|f(x)|/\alpha)  dx$.
Clearly the desired estimate (3.7) follows from (4.3), (4.6) and (4.7).

\subheading {Introducing  cancellation} As in standard Calder\'on-Zygmund theory we modify  the 
functions $f^{n,\nu}_q$  to introduce some cancellation.
Namely let $\{P_i\}_{i=1}^{M}$ be an orthonormal basis of the space of polynomials  of degree
$ \le 100 d$ on the unit cube $[-1/2,1/2]^d$ 
and for a cube $q$ with center $x_q$ and length $l(q)$ define
 the projection operator $\Pi_q$ by
$$\Pi_q[f](x)= \chi_q(x)\sum_{i=1}^M  P_i\bigl(\frac{x-x_q}{l(q)}\bigr)
 \int_q f(y) P_i\bigl(\frac{y-x_q}{l(q)}\bigr) \frac{dy}{l(q)^d}.
$$
Note that
$$|\Pi_q[h](x)|\le C \frac{1}{|q|}\int_q|h(y)|dy
\tag 4.8
$$
where $C$ is independent of
$h$ and $q$.

Let $$b^{n,\nu}_q(x) = f^{n,\nu}_q(x)-\Pi_q [f^{n,\nu}_q](x)$$
so that $b^{n,\nu}_q$ vanishes off  $q$ and 
$$\int_q b^{n,\nu}_q(x) x^{\alpha} dx =0 \qquad \text{ if } |\alpha|\le 100 d.
\tag 4.9$$

We observe that since the  $q$'s are  Whitney cubes for $\Omega$, we have
$$
\sum_{n,\nu}\big|\Pi_q f_q^{n,\nu}(x)\big|  \lc
\chi_q(x) \frac{1}{|q|}\int_q |f(x) | dx \lc \alpha ;
\tag 4.10$$
moreover by (4.8) 
$$\sum_{n,\nu}\|b^{n,\nu}_q\|_1 \lc 
\sum_{n,\nu}\|f^{n,\nu}_q\|_1 \lc
\int_q|f(x)| dx.\tag 4.11
$$

Now (4.7) will follow by Chebyshev's inequality
from
$$
\align
&\Big\| \sup_k \big| \sum_q \sum\Sb n, \nu:\\ k > k_{q}^{n,\nu} + \ka_n
\\ \nu<n^2
\endSb\Pi_q
[f^{n,\nu}_q]*d\sigma_k \big| \Big\|_2^2 \lc \alpha 
\int|f(x)| dx
\tag 4.12
\\&
\Big\| \sup_k \big| \sum_q \sum\Sb n, \nu:\\ k > k_{q}^{n,\nu} + \ka_n
\\ \nu<n^2
\endSb
b^{n,\nu}_q*d\sigma_k \big| \Big\|_2^2 \lc \alpha 
\int|f(x)| dx
\tag 4.13
\endalign
$$
since $\alpha\|f\|_1\lc \alpha^2 \int \Phi(f/|\alpha)dx$.

From (4.10) and the disjointness of the cubes $q$ we have
$\sum_{q, n,\nu}\big|\Pi_q f_q^{n,\nu}(x)\big|  \lc\alpha$
and the estimate (4.12) is immediate because of the positivity and
$L^2$ boundedness of the lacunary spherical maximal operator, and (4.10).

For the remainder of the paper we prove (4.13). We replace the $\sup$ by an $\ell^2$ norm
(in $k$) and use Minkowski's inequality to estimate
the left hand side of (4.13) by
$$
\Big[\sum_{i=1}^{N(d)}\Big(\sum_k\Big \| \sum_{q\in \fQ_i} \sum\Sb n,\nu:
\nu<n^2\\ k > k^{n,\nu}_q+\ka_n\endSb
b^{n,\nu}_q*d\sigma_k \Big\|_2^2\Big)^{1/2}\Big]^2
\lc   \sum_{i=1}^{N(d)} (I_i+II_i)
$$
where
 $$
\align
I_i&=
\sum_k \sum_{q\in\fQ_i}
\Big| \biginn{
\widetilde{d\sigma_k }*d\sigma_k *
\sum\Sb n,\nu:\\ k > k^{n,\nu}_q+\ka_n \\ \nu<n^2\endSb b^{n,\nu}_q}
{\sum\Sb n',\nu':\\ k > k^{n',\nu'}_q+\ka_{n'}\\ \nu'<{n'}^2\endSb
 b^{n',\nu'}_q} \Big|
\tag 4.14
\\
II_i &=2
\sum_k \sum\Sb q, q'\in\fQ_i\\ q \neq q'\\ l(q)\le l(q')\endSb
\Big|\sum\Sb n,\nu:\\ k > k^{n,\nu}_q+\ka_n \\ \nu<n^2\endSb
\sum\Sb n',\nu':\\ k > k^{n',\nu'}_{q'}+\ka_n' \\ \nu'<{n'}^2\endSb
\inn{
\widetilde{d\sigma_k }*d\sigma_k * b^{n,\nu}_q}
{b^{n',\nu'}_{q'}}\Big|;
\tag 4.15
\endalign
$$
and we shall consider separately the terms $I_i$ and $II_i$.

In order to estimate these expressions
we use the following well known estimate
$$
\big|\partial^{\gamma}\big[\widetilde {d\sigma_k}*d\sigma_k\big](x)\big|
\lc 2^{-k(d-1)} |x|^{-1-|\gamma|} \chi_{|x|\le 2^{k+1}}
\tag 4.16
$$
for all multiindices $\gamma\in \Bbb N_0^d$; here we need the assumption on the small support of $d\sigma_0$.

The cancellation of the functions $b^{n,\nu}_q$ will only play a role for the estimation 
of $II_i$; here no geometric information on the sets $E^{n,\nu}_q$ is used.
We carry out this estimate 
and use the moment conditions of order $N=10d$ on
the $b^{n,\nu}_q$, the fact that the cubes in $\fQ_i$ are separated   and the estimate (4.16).

Since the doubly expanded cubes are disjoint by construction of the family $\fQ_i$ we obtain 
for $l(q)<l(q')$, $x\in  q'$,
$$
|\widetilde{d\sigma_k }*d\sigma_k * b^{n,\nu}_q(x)|
\lc 2^{-k(d-1)} l(q)^{N}|x-x_q|^{-N-1} \|b^{n,\nu}_q\|_1.
$$

Thus
$$\align II_i&\lc
\sum_{n,\nu}\sum_{q\in \fQ_i} \|b^{n,\nu}_q\|_1
\sum\Sb q'\in \fQ_i\\ l(q') \ge l(q) \\q'\neq q\endSb
\sum_{k>l(q)} 2^{-(d-1)k} l(q)^{N}\dist (q,q')^{-N-1} \sum_{n',\nu'}\|b^{n',\nu'}_{q'}\|_1
\endalign$$
Next note that by (4.11)
% and the disjointness of the sets $E^{n',\nu'}_{q'}$  
$$\sum_{n',\nu'}\int_{q'} |b^{n',\nu'}_{q'}(y)| dy
\lc \int_{q'}|f(y)| dy \lc \alpha |q'|
.$$
 Moreover, we have  for  fixed $q$
that $ \dist(q,q')\gc l(q')$ if $q', q\in \fQ_i, l(q')\ge l(q)$ 
and $q\neq q'$. Thus
$$\align\sum\Sb q'\in \fQ_i\\ l(q') \ge l(q) \\ q'\neq q\endSb
&\sum_{k>l(q)} 2^{-(d-1)k} l(q)^{N}\dist (q,q')^{-N-1}|q'|
\\&\lc  l(q)^{N-d+1}\int_{|x-x_{q}|\ge l(q)} |x - x_{q}|^{-N-1} dx \lc 1.
\endalign$$
Combining the two previous estimates and applying (4.11) again  yields
$$
II_i\lc \alpha \sum_{n,\nu}\sum_q \|b^{n,\nu}_q\|_1\lc
\alpha \|f\|_1.
\tag 4.17
$$

\subheading{Estimation of the main term $\boldkey I_{\boldkey i}$}
We estimate
$I_i\le I_{i,1}+ I_{i,2}$
where
$$
\aligned
I_{i,1}&=
\sum_k \sum_{q\in\fQ_i}
\Big| \biginn{
\widetilde{d\sigma_k }*d\sigma_k *
\sum\Sb n,\nu:\\ k > k^{n,\nu}_q+\ka_n \\ \nu<n^2\endSb f^{n,\nu}_q}
{\sum\Sb n',\nu':\\ k > k^{n',\nu'}_q+\ka_{n'}\\ \nu'<{n'}^2\endSb
 b^{n',\nu'}_q} \Big|
\\
I_{i,2}&=
\sum_k \sum_{q\in\fQ_i}
\Big| \biginn{
\widetilde{d\sigma_k }*d\sigma_k *
\sum\Sb n,\nu:\\ k > k^{n,\nu}_q+\ka_n \\ \nu<n^2\endSb \Pi_q[f^{n,\nu}_q]}
{\sum\Sb n',\nu':\\ k > k^{n',\nu'}_q+\ka_{n'}\\ \nu'<{n'}^2\endSb
 b^{n',\nu'}_q} \Big|
\endaligned
\tag 4.18
$$
The estimation of $I_{i,2}$ is rather straightforward.
By (4.16) for $\gamma=0$ and by (4.10/11) we get
$$
\align
I_{i,2}
&\lc \sum_{q}
 \int_{q} \sum_{k\ge l(q)}  2^{-k(d-1)}
 \int_q \frac{1}{|x-y|}\alpha \chi_q(y) dy\,
\,\Big|\sum\Sb n',\nu':\\ k > k^{n',\nu'}_q+\ka_{n'}\\ \nu'\le {n'}^2\endSb
 b^{n',\nu'}_q(x)\Big| dx
\\
&\lc \alpha  \sum_q
\sum_{k\ge l(q)}  (2^{-k}l(q))^{d-1}
\sum\Sb n',\nu':\\ k > k^{n',\nu'}_q+\ka_{n'}\\ \nu'\le {n'}^2\endSb
\int_q|b_q^{n',\nu'}(x)|dx \lc \alpha \|f\|_1 .
\tag 4.19
\endalign
$$

Now we estimate the more substantial term  $I_{i,1}$ and use
the estimation in terms of the 
thickness of the sets
$E^{n,\nu}_q$.

Using  (4.16) with $\gamma=0$  we bound
$$I_{i,1}\lc
\sum_{k} \sum_{q\in \fQ_i }
\sum\Sb n',\nu':\\ k> k^{n',\nu'}_q+\ka_{n'}\endSb
\int |b^{n',\nu'}_q(x)|\sum_{n\ge 10} \big[A^{k,q}_{n,}(x)+
B^{k,q}_{n}(x)\big] dx
$$
where
$$\align 
A^{k,q}_{n}(x)&=
\int\limits\Sb 2^{-2n}l(q)\le \\|x-y|\le  l(q) \endSb \frac{2^{-k(d-1)}}{|x-y|}
2^n\alpha  \sum\Sb \nu<n^2:\\
k> k^{n,\nu}_q+\ka_{n}\endSb
\chi_{E^{n,\nu}_q}(y) dy,
\\
B^{k,q}_{n}(x)&=
\int\limits\Sb|x-y|\le\\ 2^{-2n} l(q)\endSb
\frac{2^{-k(d-1)}}{|x-y|}
2^n\alpha  \sum\Sb \nu<n^2:\\
k> k^{n,\nu}_q+\ka_{n}\endSb
\chi_{E^{n,\nu}_q}(y) dy .
\endalign
$$

The estimate for the terms 
involving $B^{k,q}_{n}(x)$ is straightforward. 
Since $2^{k^{n,\nu}_q}\ge l(q)$ we get 
$$
%align 
B^{k,q}_{n}(x)\lc
\int_{|x-y|\le 2^{-2n} l(q)} 2^n \alpha \frac{ 2^{-k(d-1)}}
{|x-y|}
dy 
\Big\|\sum_{\nu<n^2}
\chi_{E^{n,\nu}_q}\Big\|_\infty
\lc \alpha 
2^{-n(2d-3)} (2^{-k}l(q))^{d-1}
%\endalign
$$
and thus
$$\align
&\sum_{k} \sum_{q\in \fQ_i }
\sum\Sb n',\nu':\\ k> k^{n',\nu'}_q+\ka_{n'}\endSb
\int |b^{n',\nu'}_q(x)|\sum_{n\ge 10}  B^{k,q}_{n}(x)
dx
\\
&\lc \alpha  
 \sum_q \sum_{n', \nu'} 
\int | b^{n',\nu'}_q(x)|dx
\lc \alpha  \|f\|_1.
\endalign
$$
Next for the main term
we use that for 
$x\in q$
$$
\int\limits\Sb
2^{-n+m-1}l(q)<\\ |x-y|\le 2^{-n+m} l(q)\endSb
\frac{1}{|x-y|}
\chi_{E^{n,\nu}_q}(y) dy
\lc \Theta(E^{n,\nu}_q)
$$
and therefore
$$
A^{k,q}_{n}(x)\lc 2^{-k(d-1)} 2^n\alpha \sum\Sb
\nu\le n^2\\ k>k_q^{n,\nu}+\kappa_n\endSb 
2n\Theta(E^{n,\nu}_q) \lc
2n\alpha \sum\Sb
\nu\le n^2\\ k>k_q^{n,\nu}+\kappa_n\endSb 
2^{-(k-k_q^{n,\nu})(d-1)} .
$$
Now we perform the $k$ summation and since
$2^{-\kappa_n}=n^{-100}$ we get
$$\align
&\sum_{k} \sum_{q\in \fQ_i }
\sum\Sb n',\nu':\\ k> k^{n',\nu'}_q+\ka_{n'}\endSb
\int |b^{n',\nu'}_q(x)|\sum_{n\ge 10}  A^{k,q}_{n}(x)
dx
\\
&\lc\sum_{n'} \sum_{q\in \fQ_i }
\sum\Sb n',\nu'\endSb
\int |b^{n',\nu'}_q(x)|dx\sum_{n\ge 10} n^{3-100(d-1)}  \alpha
\lc \alpha\|f\|_1 .
\endalign
$$

Putting the estimates together we obtain
$I_i\lc \alpha\|f\|_1$ and by (4.17) the expression 
$II_i$ satisfies the same bound. This 
yields the desired estimate (4.13) and finishes the proof.
\qed

\enddocument